\documentclass{amsart}

\usepackage{amssymb}

\usepackage[all]{xy}

\newtheorem{thm}{Theorem}%[section]

\newtheorem{lem}[thm]{Lemma}

\newtheorem{prop}[thm]{Proposition}

\newtheorem{conj}[thm]{Conjecture}
\theoremstyle{definition}
\newtheorem{defn}[thm]{Definition}

\newtheorem{say}[thm]{}
\newtheorem{exmp}[thm]{Example}

\newtheorem{prob}[thm]{Problem}
\newtheorem{ques}[thm]{Question}    %!!!!!!!!!!!!!!!!!!!!

\newtheorem{rem}[thm]{Remark}          

\newtheorem*{ack}{Acknowledgments}      % \renewcommand{\theack}{} 

\newtheorem{defn-thm}[thm]{Definition--Theorem}  %!!!!!!!!!!!!!!!!!!!!!!!!
\newtheorem{defn-lem}[thm]{Definition--Lemma}  %!!!!!!!!!!!!!!!!!!!!!!!!
\newtheorem{comm}[thm]{Comment}

\theoremstyle{remark}

\renewcommand{\c}[0]{{\mathbb C}}  

\renewcommand{\o}[0]{{\mathcal O}} 
\newcommand{\z}[0]{{\mathbb Z}}

\newcommand{\p}[0]{{\mathbb P}}

\newcommand{\q}[0]{{\mathbb Q}}
\newcommand{\map}[0]{\dasharrow}
\newcommand{\qtq}[1]{\quad\mbox{#1}\quad}
\newcommand{\spec}[0]{\operatorname{Spec}}

\newcommand{\rank}[0]{\operatorname{rank}}
\newcommand{\mult}[0]{\operatorname{mult}}

\newcommand{\supp}[0]{\operatorname{Supp}}    
    
\newcommand{\codim}[0]{\operatorname{codim}}

\newcommand{\cent}[0]{\operatorname{center}}

\newcommand{\sing}[0]{\operatorname{Sing}}    
\newcommand{\ex}[0]{\operatorname{Ex}}

\newcommand{\rdown}[1]{\lfloor{#1}\rfloor}

\newcommand{\tsum}[0]{\textstyle{\sum}}

\newcommand{\dd}[0]{{\mathbb D}}
\newcommand{\bdd}[0]{\overline{\mathbb D}}

\newcommand{\sharc}[0]{\operatorname{ShArc}}

\newcommand{\farc}[0]{\widehat{\operatorname{Arc}}}

\newcommand{\ddiv}[0]{\operatorname{Div}}

\def\loccoh#1.#2.#3.#4.{H^{#1}_{#2}(#3,#4)}

\DeclareMathAlphabet{\mathchanc}{OT1}{pzc}%
                                {m}{it}

\begin{document}

\bibliographystyle{amsalpha}

%\today

\title{Arc spaces of $cA$-type singularities}
\author{Jennifer M.\ Johnson and J\'anos Koll\'ar}

\maketitle

Let $X$ be a  complex variety or an analytic space 
and  $x\in X$ a point. 
A formal arc through $x$ is a morphism $\phi:\spec \c[[t]]\to X$
such that $\phi(0)=x$.
The set of formal arcs through $x$ --
denoted by 
$\farc(x\in X)$ --
is  naturally a (non-noetherian) scheme.

A preprint of Nash, written in 1968 but only published later as
\cite{nash-arc}, describes an injection -- called the {\it Nash map} -- from the
irreducible components of $\farc(x\in X)$ to the set of
so called {\it essential divisors.}  These are the divisors
whose center on 
every resolution $\pi:X'\to X$ is an irreducible component of
$\pi^{-1}(x)$.
 The {\it Nash problem} asks if this map is 
also surjective or not. Surjectivity fails in dimensions
$\geq 3$  \cite{MR2030097, df-arc} but holds in dimension 2 \cite{boba-pp}.

In all dimensions, the most delicate cases are singularities whose
resolutions contain many rational curves. 
For example, while it is very easy to describe all arcs and their
deformations on Du Val singularities of type $A$, 
the type $E$ cases have been
notoriously hard to treat  \cite{ps-E6, pereira-phd}.

The first aim of this note is to 
determine the irreducible components of the arc space of 
$cA$-type singularities in all dimensions. 
In Section \ref{sec.irreds}  we prove the following using
 quite elementary arguments.

\begin{thm}\label{cA.sharc.thm}
Let $f(z_1,\dots, z_n)$ be a holomorphic function
whose multiplicity at the origin is  $m\geq 2$.
 Let $X:=\bigl(xy=f(z_1,\dots, z_n)\bigr)\subset \c^{n+2}$
 denote the corresponding 
 $cA$-type singularity.
 Assume that $\dim X\geq 2$.
\begin{enumerate}
\item  $\farc(0\in X)$ has $(m-1)$
irreducible  components  $\farc_i(0\in X)$ for  $0<i<m$.
\item There are dense, open subsets
 $\farc_i^{\circ}(0\in X)\subset \farc_i(0\in X)$
such that
$$
\bigl(\psi_1(t), \psi_2(t),\phi_1(t), \dots, \phi_n(t)\bigr)\in 
\farc_i^{\circ}(0\in X)
%\eqno{(\ref{cA.sharc.thm}.1)}
$$ 
iff
$\mult\psi_1(t)=i,\ \mult \psi_2(t)=m-i $ and 
$\mult f\bigl(\phi_1(t), \dots, \phi_n(t)\bigr)=m$.
\end{enumerate}
\end{thm}

We found it much harder to compute the set of essential divisors
and we have results only if $\mult_0f=2$. If $\dim X=3$ then, 
after a coordinate change, we can write the equation as
$(xy=z^2-u^m)$.
Already \cite{nash-arc} proved that these singularities
 have at most 2 essential divisors: an easy one obtained by blowing-up 
the origin and a 
 difficult one  obtained
by blowing-up the origin twice.
In Section \ref{sec.essent} 
we use ideas of \cite{df-arc} to determine the cases when
the second divisor is  essential.
The following  is obtained by  combining
 Theorem \ref{cA.sharc.thm} and Proposition \ref{res.prop}.

\begin{exmp}\label{main.nash.exmp} 
For the singularities $X_m:=(xy=z^2-u^m)\subset \c^4$
the Nash map is not  surjective  
 for  odd $m\geq 5$ but surjective for
 even $m$  and for $m=3$.
\end{exmp}

Thus the simplest counter example to the Nash conjecture is
the singularity
$$
(x^2+y^2+z^2+t^5=0)\subset \c^4.
$$
In higher dimensions our answers are less complete.
We describe  the situation for the divisors obtained
by the first and second blow-ups as above, but we do not
control other exceptional divisors. 
Using  Theorem \ref{cA.sharc.thm} 
and Proposition \ref{res.higher.prop} we get the following
partial generalization of Example \ref{main.nash.exmp}.

\begin{exmp} \label{res.higher.cor}  Let $g(u_1,\dots, u_r) $
be an analytic function near the origin. Set $m=\mult_0g$
and let $g_m$ denote the degree $m$ homogeneous part of $g$.
If $m\geq 4$ and the Nash map is  surjective  for the singularity
$$
X_g:=\bigl(xy=z^2-g(u_1,\dots, u_r)\bigr)\subset \c^{r+3}
$$
then $g_m(u_1,\dots, u_r)$  is  a perfect square.
\end{exmp}

Since we do not determine all essential divisors,
 the cases when $g_m(u_1,\dots, u_r)$ 
is  a perfect square remain undecided.

On the one hand, this can be interpreted to mean that the 
 Nash conjecture hopelessly fails in dimensions $\geq 3$.
On the other hand, the proof leads to a reformulation of
the Nash problem and to an approach that might be 
feasible, at least in dimension 3; see Section \ref{sec.revised}.

In Section \ref{sec.short} we observe that the deformations
constructed in Section \ref{sec.irreds}
also lead to an enumeration of the irreducible  components
of the space of short arcs -- introduced in \cite{k-short} --
for $cA$-type singularities.

\begin{ques}[Arcs on $cDV$ singularities]\label{cDV.ques}
It is easy to see that  Theorem \ref{cA.sharc.thm} 
is equivalent to saying that the image of every general arc
on $X$ is contained in an $A$-type surface section of $X$.

It is natural to ask if this holds for all $cDV$ singularities.
That is, let $(0\in X)\subset \c^n$ be a hypersurface singularity
such that $X\cap L^3$ is a Du~Val singularity for every general
3-dimensional linear space  (or smooth 3--fold) $0\in L^3\subset \c^n$.

Let $\phi$ be a general arc on $X$. 
Is it true that there is a 3--fold $L^3\subset \c^n$
containing the image of $\phi$ 
such that $X\cap L^3$ is a Du~Val  singularity?
\end{ques} 

 \begin{ack}
We thank  V.~Alexeev, T.~de~Fernex, R.~Lazarsfeld, 
C.~Pl\'enat and M.~Spivakovsky
for  corrections and helpful discussions. 
Partial financial support to JK
  was provided  by  the NSF under grant number 
DMS-07-58275 and by the Simons Foundation.
 Part of the paper was written while the authors visited Stanford University.
\end{ack}

\section{Arcs on $cA$-type singularities}\label{sec.irreds}

\begin{defn}[$cA$-type singularities]
In some coordinates write a hypersurface singularity as 
$$
X:=\bigl(f(x_1,\dots, x_{n+1})=0\bigr)\subset \c^{n+1}.
$$
 Assume that $X$ is singular at
the origin and 
let $f_2$ denote the quadratic part of $f$. 
If $\mult_0 f=2$ then
$(f_2=0)$ is the tangent cone of $X$ at the origin.
We say that $X$ has {\it $cA$-type} if  $\rank f_2\geq 2$
and  {\it $cA_1$-type}
that $\rank f_2\geq 3$.
By the Morse lemma, if $\rank f_2=r$ then 
we can choose local analytic or formal coordinates $y_i$ such that
$$
f=y_1^2+\cdots+y_r^2+g(y_{r+1},\dots, y_{n+1})
\qtq{where} \mult_0 g\geq 3.
$$
In the sequel we also use other forms of the quadratic part if that 
is more convenient.

Note that by adding 2 squares in new variables we get a
map
from hypersurface singularities in dimension $n-2$
(modulo isomorphism) to  $cA$-type hypersurface singularities in dimension 
$n$ (modulo isomorphism). This map is one-to-one and onto;
see \cite[Sec.11.1]{avg}. Thus $cA$-type singularities are quite
complicated in large dimensions.
\end{defn}

We rename the coordinates and write a $cA$-type singularity as
$$
X:=\bigl(xy=f(z_1,\dots, z_n)\bigr).
$$
Thus an arc through the origin is written as 
$$
t\mapsto  \bigl(\psi_1(t), \psi_2(t),\phi_1(t), \dots, \phi_n(t)\bigr)
$$
where $\psi_i, \phi_j$ are  power series such that 
$\mult \psi_i,\mult \phi_j\geq 1$ for $i=1,2$ and $j=1,\dots, n$.
We set  ${\vec \phi}(t)=\bigl(\phi_1(t), \dots, \phi_n(t)\bigr)$.

A deformation of ${\vec \phi}(t)$ is given by power series
$ \bigl(\Phi_1(t,s), \dots, \Phi_n(t,s)\bigr)$. Then we compute
$$
f\bigl(\Phi_1(t,s), \dots, \Phi_n(t,s)\bigr)\in \c[[t,s]]
$$
and try to factor it to obtain 
$$
\Psi_1(t,s) \Psi_2(t,s)=
f\bigl(\Phi_1(t,s), \dots, \Phi_n(t,s)\bigr)\in \c[[t,s]].
$$
Usually this factoring is not possible, but  Newton's method of rotating rulers
says that 
$$
f\bigl(\Phi_1(t,s^r), \dots, \Phi_n(t,s^r)\bigr) 
$$
 factors for some $r\geq 1$.

\begin{say}[Proof of Theorem \ref{cA.sharc.thm}] 
 After a linear change of coordinates we may assume that
$z_1^m$ appears in $f$ with nonzero constant coefficient.

Set $D:=\mult_t f\bigl(\phi_1(t), \dots, \phi_n(t)\bigr)$.
Assume first that $D<\infty$ and consider
$$
F(t, s):=f\bigl(\phi_1(t)+st, \phi_2(t), \dots, \phi_n(t)\bigr)=
\sum_i \frac{\partial^i f}{\partial z_1^i}\bigl(\vec\phi\bigr)
\cdot \frac{(st)^i}{i!}.
$$
We know that $t^m$  divides $F(s,t)$ (since $\mult_0 f= m$) and
$(st)^m$ appears in $F$ with nonzero coefficient
(since $z_1^m$ appears in $f$ with nonzero coefficient).
Thus $t^m$ is the largest $t$-power that divides $F(s,t)$.

Furthermore, $t^D$ is the smallest $t$-power that
 appears in $F$ with nonzero  constant coefficient. 
Thus, by Lemma \ref{newton.lem},
there is an $r\geq 1$ such that
$$
F(t, s^r)=u(t,s)\prod_{i=1}^D \bigl(t-\sigma_i(s)\bigr)
$$
where $u(0,0)\neq 0$ and  $\sigma_i(0)=0$. Furthermore, 
 exactly $m$ of the $\sigma_i$ are identically zero.

For $j=1,2$ write  $\psi_j(t)=t^{a_j}v_j(t)$
where $v_j(0)\neq 0$. Note that $a_1+a_2=D$
and $u(t,0)=v_1(t)v_2(t)$.

Divide $\{1,\dots, D\}$ into two disjoint subsets $A_1, A_2$
such that $|A_j|=a_j$ and they both contain at least 1 index
$i$ such that $\sigma_i(t)\equiv 0$.
Finally set
$$
\Psi_1(t,s)=v_1(t)\cdot \prod_{i\in A_1} \bigl(t-\sigma_i(s)\bigr)
\qtq{and}
\Psi_2(t,s)=\frac{u(t,s)}{v_1(t)}
\cdot \prod_{i\in A_2} \bigl(t-\sigma_i(s)\bigr).
$$
Then
$$
\bigl(\Psi_1(t,s), \Psi_2(t,s), \phi_1(t)+st, \phi_2(t), \dots, \phi_n(t)\bigr)
$$
is a deformation of 
$\bigl(\psi_1(t), \psi_2(t),\phi_1(t), \dots, \phi_n(t)\bigr)$
whose general member is in the $r$th 
irreducible  component as in (\ref{cA.sharc.thm}.2)
iff exactly $r$ of the $\{\sigma_i: i\in A_1\}$ are identically zero.

(This also shows that arcs with  $\mult \psi_1(t)\geq m-1$ and
$\mult \psi_2(t)\geq m-1$
constitute  the intersection of all of the irreducible  components.)

If $D=\infty$, that is, when $f\bigl(\phi_1(t), \dots, \phi_n(t)\bigr)$
is identically zero, we need to perform some
similar preliminary deformations first.

First, if both $\psi_1(t), \psi_2(t)$ are identically zero
then we can take
$$
 \bigl(st, 0, \phi_1(t), \phi_2(t), \dots, \phi_n(t)\bigr).
$$
Hence, up-to interchanging $x$ and $y$, we may assume  that
$d:=\mult \psi_1(t)<\infty$. 
Again assuming that $z_1^m$ appears in $f$ with nonzero coefficient,
we see that
$$
F(t,s):=f\bigl(\phi_1(t)+st^{d+1}, \phi_2(t), \dots, \phi_n(t)\bigr)
$$
is not identically zero and divisible by $t^{d+1} $.
Thus $F(t,s)/\psi_1(t) $  is holomorphic and divisible by $t$. Therefore
$$
\Bigl(\psi_1(t), \frac{F(t,s)}{\psi_1(t)}, 
\phi_1(t)+st^{d+1}, \phi_2(t), \dots, \phi_n(t)\Bigr)
$$
is a deformation of 
$\bigl(\psi_1(t), 0, \phi_1(t), \phi_2(t), \dots, \phi_n(t)\bigr)$
such that
$$
\mult_t f\bigl(\phi_1(t)+st^{d+1}, \phi_2(t), \dots, \phi_n(t)\bigr)<\infty
$$
for $0<|s|\ll 1$.
\qed
\end{say}

We used Newton's lemma on Puiseux series solutions in the following form.

\begin{lem}\label{newton.lem}
Let $g(x,y)\in \c[[x,y]]$ be a power series.
Assume that $m:=\mult_0g(x,0)<\infty$. Then there is an $r\geq 1$
 such that one can write $g(x, z^r)$ as
$$
g(x, z^r)=u(x,z)\prod_{i=1}^m \bigl(x-\sigma_i(z)\bigr)
$$
where $u(0,0)\neq 0$ and $\sigma_i(0)=0 $
for every $i$.
The representation is unique, up-to permuting the $\sigma_i(z) $.

Furthermore, if 
 $g(x,y)$ is holomorphic on the bidisc $\bdd_x\times \dd_y$
then $u(x,z)$ and the $\sigma_i(z)$ 
are holomorphic on the smaller bidisc 
$\bdd_x\times \dd_z(\epsilon)$ for some $0<\epsilon\leq 1$.
\qed
\end{lem}

\section{Essential divisors on $cA_1$-type 3-fold singularities}\label{sec.essent}

In dimension 3, the only  $cA_1$-type singularities are
 $X_m:=(xy=z^2-t^m)$ for $m\geq 2$.
Already \cite[p.37]{nash-arc} proved that they have
at most 2 essential divisors. We use the method of \cite[4.1]{df-arc}
to determine the precise count.

\begin{defn} \label{first.ess.defn}
Let $X$ be a normal variety or analytic space
and $E$ a divisor over $X$. That is, there is a 
 birational or bimeromorphic morphisms
$p:X'\to X$ such that $E\subset X'$ is an exceptional divisor.
The closure of $p(E)\subset X$ is called the {\it center}
of $E$ on $X$; it is denoted by $\cent_XE$.
If $\cent_XE=\{x\}$, we say that $E$ is a divisor over $(x\in X)$. 

We say that $E$ is an {\it essential divisor} over $X$ if
for every resolution of singularities $\pi:Y\to X$,
 $\cent_YE$  is an irreducible component
of $\pi^{-1}\bigl(\cent_XE\bigr)$.
(Note that $\pi^{-1}\circ p:X'\map Y$ is regular on a dense subset of $E$,
hence $\cent_YE$ is defined.)

If $X$ is an analytic space, then $Y$ is allowed to be any
analytic resolution. If $X$ is algebraic, one gets
slightly different notions depending on whether one 
allows $Y$ to be a quasi-projective variety, an algebraic space
or an  analytic  space;
 see \cite{df-arc}.
We believe that for the Nash problem it is natural to allow
 analytic resolutions.
\end{defn}

\begin{prop} \label{res.prop}
Set $X_m:=(xy=z^2-t^m)\subset \c^4$. 
\begin{enumerate}
\item If  $m\geq 5$ is odd, there are 2 essential divisors.
\item If $m\geq 2$ is even or $m= 3$, there is 1 essential divisor.
\end{enumerate}
\end{prop}

Even in dimension 3, it seems surprisingly difficult
to determine the set of essential divisors.
A basic invariant is given by the discrepancy.

\begin{defn} \label{discrep.defn}
Let $X$ be a normal variety or analytic space.
Assume for simplicity that the canonical class $K_X$ is Cartier.
(This holds for all hypersurface singularities.)
Let $\pi:Y\to X$ be a 
resolution of singularities  and write
$$
K_Y\sim \pi^*K_X+\tsum_i a(E_i,X)E_i
$$
where the $E_i$ are the $\pi$-exceptional divisors.
The integer $a(E_i,X)$ is called the 
 {\it discrepancy}   of  $E_i$.
(See \cite[Sec.2.3]{km-book} for basic references and more general definitions.)

For example, let $X$ be smooth and $Z\subset X$ a smooth subvariety of 
codimension $r$. Let $\pi_Z:B_ZX\to X$ denote the blow-up
and $E_Z\subset B_ZX$ the exceptional divisor.
Then $a(E_Z, X)=r-1$ and easy induction shows that
$a(F,X)\geq r$ for every other divisor whose center on $X$ is $Z$.

We say that $X$ is {\it canonical} (resp.\ {\it terminal})
of  $a(E_i,X)\geq 0$ (resp.\ $a(E_i,X)> 0$) for every
resolution and every exceptional divisor. 

For instance, normal $cA$-type singularities are canonical
and a $cA$-type singularity is terminal iff
its singular set has codimension $\geq 3$; see
\cite{Reid83} for a proof that applies to all $cDV$ singularities or
 \cite[1.42]{kk-singbook} for a simpler argument in the $cA$ case.
\end{defn}

\begin{say}[Resolving $X_m$] 
Blow up the origin to get $\pi_1:X_{m,1}:=B_0X_m\to X_m$.
The exceptional divisor is the singular quadric
$E_1\cong (xy-z^2=0)\subset \p^3(x,y,z,t)$. $B_0X_m$ has one singular point,
 visible in the chart
$$
(x_1, y_1, z_1, t):=\bigl(x/t, y/t, z/t, t\bigr)
$$
where the local equation is $x_1y_1=z_1^2-t^{m-2}$.
We can thus blow up the origin again and continue.
After $r:=\rdown{\tfrac{m}{2}}$ steps we have a resolution
$$
\Pi_r: X_{m,r}\to X_{m,r-1}\to \cdots \to X_{m,1} \to X_m.
$$
We get $r$ exceptional divisors
$E_r,\dots, E_1$. For $1\leq c\leq r$ the divisor $E_c$
first appears on $X_{m,c}$. At the unique singular point
one can write the  local equation as
$$
X_{m,c}=\bigl(x_cy_c=z_c^2-t^{m-2c}\bigr)
\qtq{and}
E_c=(t=0).
$$
where
$(x_c, y_c, z_c, t):=\bigl(x/t^c, y/t^c, z/t^c, t\bigr)$.
\end{say}

We thus need to decide which of the divisors $E_1,\dots, E_{\rdown{\tfrac{m}{2}}}$ 
are essential. It is easy to see that $E_1$ is  essential and a
direct computation (\ref{lem.res4})
shows that  $E_3,\dots, E_{\rdown{\tfrac{m}{2}}}$
are not. 
(This is actually not needed in order to establish
 Example \ref{main.nash.exmp}.)
The hardest is to decide what happens with $E_2$.

\begin{lem}\label{lem.res1} Notation as above. Then
\begin{enumerate}
\item $a(E_c, X_m)=c$ for every $c$.
\item $E_1$ is the only exceptional divisor whose center is
the origin and whose discrepancy is $1$.
\item $E_1$ appears on every resolution of $X_m$
whose exceptional set is a divisor.
\item Let $p:Y\map X_m$ be any (not necessarily proper)
 bimeromorphic map from a smooth analytic space $Y$ such that
$\cent_YE_1\subset Y$ is not empty. Then $\cent_YE_1$ is an
irreducible component of the exceptional set $\ex(p)$. 
\end{enumerate}
\end{lem}

Proof. The first claim follows from the formula
$$
\Pi_r^*\Bigl(\tfrac{dx\wedge dy\wedge dt}{z}\Bigr)=
t^{-c}\cdot \tfrac{dx_c\wedge dy_c\wedge dt}{z_c}.
$$
Let $F$ be any other exceptional divisor whose center is
the origin. Then $\cent_{X_r}F$ lies on one of the $E_c$,
thus $a(F,X)>a(E_c,X)\geq 1$. (This also proves that $X_m$ is terminal.)

To see (3) set $W_1:=\cent_YE_1\subset Y$.
Let $F_i\subset Y$ be the exceptional divisors and note that,
as in \cite[2.29]{km-book}, 
$$
a(E_1, X_m)\geq 
\bigl(\codim_YW_1-1\bigr)+\tsum_i \mult_{W_1}F_i\cdot a(F_i, X_m).
\eqno{(\ref{lem.res1}.5)}
$$ 
Note that  $a(E_1, X_m)=1$ and
$a(F_i, X_m)\geq 1$ for every $i$. 
If $W_1$ is not  an
irreducible component of $\ex(p)$ then
$W_1\subset F_i$ form some $i$ and then both terms on 
the right hand side of (\ref{lem.res1}.5) are positive, a contradiction.
\qed

\begin{lem}\label{lem.res2} If $m\in \{2,3\}$ then $B_0X$ is smooth,
hence the only essential divisor is $E_1$. \qed
\end{lem}

\begin{say}[Small resolutions and  factoriality  of $X_m$] 
\label{small.res.say.1}
 If $m=2a$ is even, then
$X_m$ has a small resolution obtained by blowing up either
$(x=z-t^a=0)$ or $(x=z+t^a=0)$. The resulting blow-ups 
$ Y^{\pm}_{2a}\subset  \c^4_{xyzt}\times \p^1_{uv}$ are 
defined by the equations
$$
 Y^{\pm}_{2a}:=\rank \left(
\begin{array}{ccc}
x & z\pm t^a & u\\
z\mp t^a & y & v
\end{array}
\right)
\leq 1
\eqno{(\ref{small.res.say.1}.1)}
$$
By contrast,  $X_m$
does not have small resolutions if $m$ is odd.
More generally, let 
$$
X_f:=\bigl(xy=f(z,t)\bigr)\subset \c^{4}
$$
be an isolated   $cA$-type singularity.
Write $f=\prod_j f_j$ as a product of
 irreducibles. The $f_j$ are distinct since the singularity is isolated.
Set $D_j:=(x=f_j=0)$. 
By  \cite[2.2.7]{k-etc} the local divisor class group is
$$
\ddiv\bigl(0\in X_f\bigr)=
\bigl(\tsum_j \z[D_j]\bigr)\big/\tsum_j [D_j].
\eqno{(\ref{small.res.say.1}.2)}
$$
In particular,  $X_f$
is factorial  iff
$f$ is irreducible.

This formula works both algebraically and analytically.
If we are interested in the affine variety $X_f$, then
we consider factorizations of $f$ in the polynomial ring.
If  we are interested in the complex analytic germ $X_f$,
then we consider factorizations of $f$ in the ring of germs of
analytic functions. Thus, for example, 
$$
(xy=z^2-t^2-t^3)\subset \c^4
$$
is algebraically factorial,
since $z^2-t^2-t^3$ is an irreducible polynomial,
but it is not analytically factorial,
since 
$$
z^2-t^2-t^3=\bigl(z-t\sqrt{1+t}\bigr)\bigl(z+t\sqrt{1+t}\bigr).
$$

Thus if $m$ is odd then $X_m$ is factorial  
(both algebraically and analytically) and it
does not have small resolutions; see Lemma \ref{purecodim1.lem}
for stronger results.
\end{say}

\begin{lem}\label{lem.res3}  If $m$ is even then 
there is a divisorial resolution whose sole exceptional divisor is
 birational to $E_1$.
Thus the only essential divisor is $E_1$. 
\end{lem}

Proof.  The $m=2$ case is in (\ref{lem.res2}), hence we may assume that
$m=2a\geq 4$. 

There are 2 ways to obtain such resolutions.
First, we can blow up the exceptional curve
in  either of the $Y^{\pm}_{2a} $ as in (\ref{small.res.say.1}.1).

Alternatively, we
 first blow up the origin to get  $B_0X_m$ which has one singular point 
with local
equation $x_1y_1=z_1^2-t_1^{2a-2}$
and then   blow up
$D^+:=(x_1=z_1+t_1^{a-1}=0)$ or  $D^-:=(x_1=z_1-t_1^{a-1}=0)$.\qed

\begin{lem}\label{lem.res4} \cite[p.37]{nash-arc}
The divisors $E_3,\dots, E_r$ are not essential. 
\end{lem}

Proof. If $m$ is even, this follows from
(\ref{lem.res3}), but for the proof below the parity of $m$
does not matter.

If $2b\geq a\geq 0$ and $m\geq a$ then 
  $(u,v,w,t)\mapsto  (ut, vt^{a+1}, wt^{b+1}, t)=(x,y,z,t)$ 
defines a birational map
$$
g(a,b,m):Z_{abm}:=(uv=w^2t^{2b-a}-t^{m-2-a}\bigr)\to 
X_m.
$$
Note that 
$\ex\bigl(g(a,b,m)\bigr)=(t=0)$ is mapped to the origin and 
$Z_{abm}$ is smooth along the $v$-axis, save at the origin.

If $1\leq c\leq m/2$ then
  $(x_c,y_c,z_c,t)\mapsto  (x_ct^c, y_ct^c, z_ct^c, t)=(x,y,z,t)$ 
defines  a birational map
$$
h(c,m): X_{m,c}:=(x_cy_c=z_c^2-t^{m-2c}\bigr)\to 
X_m.
$$
By composing we get a birational map  
$g(a,b,m)^{-1}\circ h(c,m): Y_c\map Z_{abm}$
given by
$$
(x_c,y_c,z_c,t)\mapsto 
(x_ct^{c-1}, y_ct^{c-a-1}, z_ct^{c-b-1}, t)= (u,v,w,t)
$$
which is a morphism if $c\geq a+1,b+1$. 
If $c=a+1$ and $c>b+1$ then we have
$$
(x_c,y_c,z_c,t)\mapsto 
(x_ct^{c-1}, y_c, z_ct^{c-b-1}, t)= (u,v,w,t)
$$
which maps  $E_c$ to the $v$-axis.

If $c\geq 3$ then by setting $a=c-1, b=c-2$ we get
a birational  morphism 
$p(c,m):=g(c,c{-}1,m)^{-1}\circ h(c,m)$
given by 
$$
(x_c,y_c,z_c,t)\mapsto 
(x_ct^c, y_c, z_ct, t)= (u,v,w,t).
$$
Note that
$$
p(c,m): 
 Y_{c}=(x_cy_c=z_c^2-t^{m-2c}\bigr)\to 
(uv=w^2t^{c-2}-t^{m-c}\bigr)=Z_{c,c-1,m}
$$
 maps  $E_c$ onto the $v$-axis.
Thus $E_c$ is not essential for $c\geq 3$. \qed

\begin{lem}\label{lem.res5}  
If $m\geq 5$ is odd then $E_2$ is essential.
\end{lem}

Proof. 
We follows the arguments in  \cite[4.1]{df-arc}.
 Let $p:Y\to X_m$ be any resolution 
and set $Z:=\cent_YE_2\subset Y$. 
Since $X_m$ is factorial (here we use that $m$ is odd),
 $\ex(p)$ has pure dimension 2 by (\ref{purecodim1.lem}.2).

Assume to the contrary that $Z$ is not a divisor.
Using that $a(E_2, X_m)=2$, (\ref{lem.res1}.5)
implies that  $Z$ is a curve,
 there is a unique exceptional divisor $F\subset Y$ that contains
$Z$, $F$ is  smooth at  general points of $Z$ and $a(F,X_m)=1$.

If $p(F)$ is a curve then $Z$ is an irreducible component
of $p^{-1}(0)$.
The remaining case is when $p(F)=0$, thus $F=E_1$ by (\ref{lem.res1}.2).

Since $t$ vanishes along $E_2$ with multiplicity 1,
it also vanishes along $Z$ with multiplicity 1.
Since $p^*x, p^*y, p^*z, p^*t$ all vanish along $E_1$ 
the rational functions
$p^*(x/t), p^*(y/t), p^*(z/t)$
are regular  generically along $Z$. 
Thus  $p_1:=\pi_1^{-1}\circ p:Y\map X_{m,1}$
is a morphism  generically along $Z$.
Note that our $E_2$ is what we would call $E_1$ if we started
with $ X_{m,1}$. Applying (\ref{lem.res1}.4) to
$p_1:Y\map  X_{m,1}$ we see that $Z$ is an irreducible component
of $\ex(p_1)$. Since $m$ is odd, $X_{m,1}$ is
analytically factorial by (\ref{small.res.say.1}), hence
$Z$ is a divisor by (\ref{purecodim1.lem}.2). This is a contradiction.
\qed

\begin{lem} \label{purecodim1.lem}
Let $X, Y$ be  normal varieties or analytic spaces
and $g:Y\to X$ a  birational or bimeromorphic morphism.
Then the exceptional set $\ex(g)$ has pure codimension 1
in $Y$ in the following cases.
\begin{enumerate}
\item $Y$ is an algebraic variety and $X$ is $\q$-factorial.
\item $\dim Y=3$ and  $X$ is analytically locally $\q$-factorial.
\end{enumerate}
\end{lem}

Proof. The algebraic case is well known; see for instance
the method of \cite[Sec.II.4.4]{shaf}.

If $\dim Y=3$ and $\ex(g)$ does not have pure codimension 1
then it has a 1-dimensional irreducible component
$C\subset Y$. After replacing $X$ by a suitable neighborhood of
$g(C)\in X$ we may assume that there is a divisor
$D_Y\subset Y$ such that $\ex(g)\cap D_Y$ is a single point of
$C$ and $g|_{D_Y}$ is proper. Thus $D_X:=g(D_Y)$ is a divisor on $X$.
If $mD_X$ is Cartier then so is $g^*(mD_X)$ hence its support
has pure codimension 1 in $Y$. On the other hand,
$\supp\bigl(g^*(mD_X)\bigr)=\ex(g)\cup D_Y$ does not have
pure codimension 1. 
(Note that there are many possible choices for $D_Y$;
the resulting $D_X$ determine an algebraic equivalence class of
divisors.)\qed
\medskip

Somewhat surprisingly, the analog of 
(\ref{purecodim1.lem}.2) fails in dimension 4.

\begin{exmp} Let $W\subset \p^4$ be a smooth quintic 3--fold
and $C\subset W$ a line whose normal bundle is $\o(-1)+\o(-1)$.
Let $X\subset \c^5$ denote the cone over $W$ with vertex $0$;
it is analytically locally factorial by \cite[XI.3.14]{sga2}.

The exceptional divisor of the blow-up $B_0X\to X$ can be identified with
$W$; let $C\subset B_0X$ be our line. Its 
normal bundle is $\o(-1)+\o(-1)+\o(-1)$.

Blow up the line $C$ to obtain $B_CB_0X\to B_0X$.
Its exceptional divisor is $E\cong \p^1\times \p^2$. 
One can contract $E$ in the other direction to obtain
$g:Y\to X$. 

By construction, $\ex(g)$ is the union of $\p^2$ and of a
3-fold obtained from $W$ by flopping the line $C$.
The two components intersect along a line.
\end{exmp}

\begin{rem} \label{corves.to.locdiv}  We will need to understand in
detail the  proof of (\ref{purecodim1.lem}.2)
for
$$
X_c:=\bigl(xy=z^2-ct^{m}\bigr)\subset \c^4\qtq{where $c\neq 0$.}
$$
Let $g_c:Y_c\to X_c$ be a proper birational or bimeromorphic morphism
and $E_c\subset \ex(g_c)$ a 1-dimensional irreducible component.

The proof of (\ref{purecodim1.lem}.2) associates to $E_c$ an
algebraic equivalence class of non-Cartier divisors on $X_c$.
Thus $m$ has to be even by (\ref{small.res.say.1}).

If $m=2a$ is even then the divisor class group is
$\ddiv(X_c)\cong \z$. The two possible generators
correspond to
$(x=z-\sqrt{c}t^a=0)$ and  $(x=z+\sqrt{c}t^a=0)$.
Starting with $E_c$ we constructed a divisor
$D_c\subset X_c$ which is a nontrivial element of
$\ddiv(X_c)$. Thus $[D_c]$ is a positive multiple of either
$(x=z-\sqrt{c}t^a=0)$ or $(x=z+\sqrt{c}t^a=0)$.
Hence, to $E_c\subset Y_c$ we can associate a choice of
$\sqrt{c} $.

This may not be very interesting for a fixed value of $c$
(since many other choices are involved)
but it  turns out to be quite useful when $c$ varies.
\end{rem}

\begin{prop}\label{div.varies.infams}
Let $g(u_1,\dots, u_r, v) $ be a holomorphic function for
$u_i\in \c$ and $|v|<\epsilon$ such that 
 $g(u_1,\dots, u_r, 0)$ is not identically zero. Set 
$$
X:=\bigl(xy=z^2-v^mg(u_1,\dots, u_r, v)\bigr)\subset \c^{r+4}
$$
Let  $\pi:Y\to X$ be a  birational or bimeromorphic morphism.
Assume that there is an irreducible component
$Z\subset \ex(\pi)$ that dominates  $(x=y=z=v=0)\subset X$,
 has codimension  $\geq 2$ in $Y$ and such that
$\pi|_{Z}:Z\to (x=y=z=v=0)$ has connected fibers.

Then $m$ is even and $g(u_1,\dots, u_r,0)$ is a perfect square.
\end{prop}

Proof. 
For general ${\mathbf c}=(c_1,\dots, c_r)\in \c^r $ the repeated
hyperplane section
$$
X({\mathbf c}):=\bigl(xy=z^2-v^mg({\mathbf c}, v)\bigr)\subset \c^{4}
$$
has an isolated singularity at the origin and we get a
proper birational or bimeromorphic morphism
$$
\pi({\mathbf c}):Y({\mathbf c})\to X({\mathbf c})
$$
where $Y({\mathbf c})\subset Y$ is the 
 preimage of $X({\mathbf c}) $.

Furthermore,
 $Z({\mathbf c}):=Z\cap Y({\mathbf c})$ is an 
 irreducible component of $ \ex\bigl(\pi({\mathbf c})\bigr)$
and has codimension $\geq 2$  in $Y({\mathbf c})$.

Thus, as we noted above,
 $m=2a$ is even and our construction gives a function
$$
(c_1,\dots, c_r)\mapsto \mbox{a choice of } \sqrt{g(c_1,\dots, c_r,0)}.
$$
It is clear  that this function is
continuous on a Zariski open set  $U\subset \c^r$.
Therefore 
$g(u_1,\dots, u_r,0)$ is a perfect square. \qed

\begin{rem}\label{converse.to.square}
Conversely, assume that $m$ is even and $g(u_1,\dots, u_r, 0)=
h^2(u_1,\dots, u_r)$ is a square. Write
the equation of $X$ as
$$
xy=z^2-v^m\bigl(h^2(u_1,\dots, u_r)+vR(u_1,\dots, u_r, v)\bigr)
$$
Over the open set $X^0\subset X$ where $h\neq 0$, change coordinates to
$w:=h^{-2}v$. (Equivalently, blow up $(v=h=0)$ twice.) Then
$$
D:=\bigl(x=z-w^{m/2}h^{m+1}\sqrt{1+wR(u_1,\dots, u_r, h^2w)}\bigr)
$$
is a globally well defined analytic divisor.
Blowing it up
gives a bimeromorphic morphism $X_D\to X$ whose exceptional set
over $X^0$
has codimension 2.

It seems that  even if $X$ is algebraic, 
usually $X_D$ is not an algebraic variety.

\end{rem}

\section{Essential divisors on $cA_1$-type singularities}\label{sec.essent.higher}

In higher dimensions    $cA_1$-type singularities 
are more complicated and their resolutions are much harder to understand.
There is no simple complete answer as in dimension 3.

In the previous Section, the key part is to understand the
exceptional divisors that correspond to the first 2 blow-ups.
These are the 2 divisors that we understand in higher dimensions as well.

\begin{say}[Defining $E_1$ nd $E_2$]\label{def.E1.E2}
In order to fix notation, write the equation  as
$$
X:=\bigl(xy=z^2-g(u_1,\dots, u_r)\bigr)\subset \c^{r+3}.
\eqno{(\ref{def.E1.E2}.1)}
$$
Set $m:=\mult_0 g$ and let $g_s(u_1,\dots, u_r) $ denote the
homogeneous degree $s$ part of $g$.
In a typical local chart the 1st blow-up  $\sigma_1:X_1:=B_0X\to X$ is given by
$$
x_1y_1=z_1^2-\bigl(u'_r\bigr)^{-2}
g(u'_1u'_r,\dots, u'_{r-1}u'_r, u'_r)
\eqno{(\ref{def.E1.E2}.2)}
$$
where $x=x_1u'_r, y=y_1u'_r, z=z_1u'_r$,
$u_1=u'_1u'_r, \dots, u_{r-1}=u'_{r-1}u'_r$ and $ u_r=u'_r$.
The exceptional divisor is the rank 3  quadric 
$$
E_1:=\bigl(x_1y_1-z_1^2=0\bigr)\subset \p^{r+2}.
\eqno{(\ref{def.E1.E2}.3)}
$$
Note also that
$$
\begin{array}{l}
\bigl(u'_r\bigr)^{-2}g(u'_1u'_r,\dots, u'_{r-1}u'_r, u'_r)=\\
\qquad =\bigl(u'_r\bigr)^{m-2}\Bigl(g_m(u'_1,\dots, u'_{r-1}, 1)+
u'_rg_{m+1}(u'_1,\dots, u'_{r-1}, 1)+\cdots\Bigr).
\end{array}
\eqno{(\ref{def.E1.E2}.4)}
$$
From this we see that, for $m\geq 4$, the blow-up
 $X_1$ is singular along the closure of the
 linear space
$$
L:=(x_1=y_1=z_1=u'_r=0),
\eqno{(\ref{def.E1.E2}.5)}
$$
$X_1$ has terminal singularities 
and a general 3-fold section has equation
$$
x_1y_1=z_1^2-\bigl(u'_r\bigr)^{m-2}
\Bigl(g_m(c_1,\dots, c_{r-1}, 1)+
u'_rg_{m+1}(c_1,\dots, c_{r-1}, 1)+\cdots\Bigr).
$$

Blowing up the closure of $L$ we obtain $X_2$ with exceptional divisor
$E_2$. As in Lemma \ref{lem.res1} we compute that

\begin{enumerate}\setcounter{enumi}{5}
\item $a(E_1, X)=r$,
\item $a(E_2, X)=r+1$,
\item $a(F, X)\geq r+1$ for every other exceptional divisor whose center
on $X$ is the origin and
\item the pull-backs of the $u_i$ vanish along $E_1, E_2$ with
multiplicity 1.
\end{enumerate}

\end{say}

The key computation  is the following.

\begin{prop} \label{res.higher.prop} Notation as above and 
assume that $m\geq 4$.
\begin{enumerate}
\item  $E_1$ is an essential divisor.
\item $E_2$ is an essential divisor iff
 $g_m(u_1,\dots, u_r)$
is not a perfect square.
\end{enumerate}
\end{prop}

%\begin{say}[Proof of Proposition \ref{res.higher.prop}]

Proof. By (\ref{def.E1.E2}.6) and (\ref{def.E1.E2}.8),
$E_1$ has the smallest 
 discrepancy among all divisors over $X$ whose center on $X$ is
the origin. Thus $E_1$ is  essential by Proposition \ref{can.mind.ess.prop}.

If  $E_2$ is non-essential then there is a resolution
$\pi:Y\to X$ and an irreducible component
$W\subset \supp\pi^{-1}(0)$ such that
$Z:=\cent_YE_2\subsetneq W$. 
By (\ref{def.E1.E2}.9), the $\pi^*u_i$ vanish at a general point
of $Z$ with multiplicity 1. Since the $\pi^*u_i$ vanish along $W$,
this implies that $\supp\pi^{-1}(0)$ 
 is smooth at a general point of
$Z$. In particular, 
$W$ is the only irreducible component of $\supp\pi^{-1}(0)$  that
contains $Z$ and $W$  is smooth at general points of
$Z$.
Therefore the blow-up $B_WY$ is smooth
over the generic point of $Z$. So, if
we replace $Y$ by a suitable desingularization of
 $B_WY$, we get a situation as before where
in addition $W$ is a divisor. 
%In order to simplify notation, we  assume from now on that  $W$.

The $\pi^*u_i$ are local equations of $W$ at general points of $Z$
and  $\pi^*x, \pi^*y, \pi^*z$ all vanish along $W$. 
Thus the rational functions
$$
\pi^*(x/u_r), \pi^*(y/u_r),\pi^*(z/u_r), \pi^*(u_1/u_r),\dots, 
\pi^*(u_{r-1}/u_r),
$$
are all regular at general points of $Z$.
Hence the birational map  $\sigma_1^{-1}\circ \pi:Y\to B_0X=X_1$
is a morphism at general points of $Z$. Furthermore,
 $\sigma_1^{-1}\circ \pi $ maps $W$ birationally to $E_1\subset X_1$
and it is not a local isomorphism along $Z$ since $Y$ is smooth
but $X_1$ is singular along the center $L$ of $E_2$.
Thus $Z$ is an irreducible component of 
$\ex\bigl(\sigma_1^{-1}\circ \pi\bigr)$. 
Since $E_2\to L$  has connected fibers, all the assumptions of
 Proposition \ref{div.varies.infams} are satisfied by 
the equation of the blow-up 
$$
x_1y_1=z_1^2-\bigl(u'_r\bigr)^{m-2}
\Bigl(g_m(u'_1,\dots, u'_{r-1}, 1)+
u'_rg_{m+1}(u'_1,\dots, u'_{r-1}, 1)+\cdots\Bigr).
\eqno{(\ref{res.higher.prop}.3)}
$$
Thus $m$ is even and $g_m(u'_1,\dots, u'_{r-1}, 1) $
is a perfect square. Since it is a
dehomogenization of $g_m(u_1,\dots, u_{r-1}, u_r) $,
the latter is also a  perfect square.

The converse follows from Remark \ref{converse.to.square}.
\qed

\medskip

\begin{defn} For $(x\in X)$ let $\operatorname{min-discrep}(x\in X)$
be the infimum of the discrepancies $a(E,X)$ where $E$
runs through all divisors over $X$ such that $\cent_XE=\{x\}$.
(It is easy to see that either $\operatorname{min-discrep}(x\in X)\geq -1$
and the infimum is a minimum or 
$\operatorname{min-discrep}(x\in X)=-\infty$; cf.\ \cite[2.31]{km-book}. 
We do not need these
facts.)
\end{defn}

\begin{prop} \label{can.mind.ess.prop}
Let $(x\in X)$ be a canonical singularity
and $E$ a  divisor over $X$ such that
\begin{enumerate}
\item $\cent_XE=\{x\}$ and
\item $a(E,X)<1+\operatorname{min-discrep}(x\in X)$.
\end{enumerate}
Then $E$ is essential.
\end{prop}

Proof. Let $F$ be any non-essential divisor over $X$ whose center
on $X$ is the origin. Thus there is a resolution
$\pi:Y\to X$ and an irreducible component
$W\subset \supp\pi^{-1}(x)$ such that
$$
Z:=\cent_YF\subsetneq W.
$$
Let $E_W$ be the divisor obtained by blowing up $W\subset Y$.
As we noted in (\ref{discrep.defn}), 
$$
a(E_W, Y)=\codim_YW -1\qtq{and} a(F, Y)\geq\codim_YZ -1\geq \codim_YW.
\eqno{(\ref{can.mind.ess.prop}.3)}
$$
Write  $K_Y=\pi^*K_X+D_Y$ where  $D_Y$ is effective since
$X$ is canonical and 
note  that
$$
a(E_W, X)=a(E_W, Y)+\mult_WD_Y \qtq{and} 
a(F, X)\geq a(F, Y)+\mult_ZD_Y.
\eqno{(\ref{can.mind.ess.prop}.4)}
$$
Since $\mult_ZD_Y\geq \mult_WD_Y$, we conclude that
$$
a(F, X)\geq 1+a(E_W, X)\geq 1+\operatorname{min-discrep}(x\in X).
\eqno{(\ref{can.mind.ess.prop}.5)}
$$
Thus any divisor $E$ with $a(E,X)<1+\operatorname{min-discrep}(x\in X)$
is essential. \qed
\medskip

\section{Short arcs}\label{sec.short}

Let $\dd\subset \c$ denote the open unit disk and $\bdd\subset \c$ its closure.
The open (resp.\ closed) disc of radius $\epsilon$ is denoted by
$\dd(\epsilon)$ (resp.\ $\bdd(\epsilon)$).
If several variables are involved, we use a subscript to indicate
the name of the coordinate.

\begin{say}[Short arcs]\cite{k-short}\label{short.defs.say}
  Let $X$ be an analytic space and $p\in X$ a  point.
A {\it short arc} on $(p\in X)$ is a 
holomorphic map $\phi(t): \bdd_t\to X$ such that
$\supp\phi^{-1}(p)=\{0\}$.

The space of all short arcs is denoted by $\sharc(p\in X)$.
It has a natural topology and most likely also a
complex structure that, at least for isolated singularities,
 locally can be written as the product of 
a (finite dimensional) complex space and of a complex Banach space;
see \cite[Sec.11]{k-short} for details.

A {\it  deformation} of  short arcs  is a
holomorphic map $\Phi(t,s): \bdd_t\times \dd_s\to X$ such that
$\Phi(t,s_0): \bdd_t\to X$ is a short arc for every $s_0\in \dd_s$.
Equivalently, if
$\supp\Phi^{-1}(p)=\{0\}\times \dd_s$.

In general the space of short arcs has
more connected components 
than the space of formal arcs. As a simple example,
consider arcs on $(xy=z^n)\subset \c^3$.   For $0<i<m$
the  deformations
$$
(t,s)\mapsto  \bigl(t^i(t+s)^{m-i}, t^{m-i}(t+s)^{i}, t(t+s)\bigr)
\eqno{(\ref{short.defs.say}.1)}
$$
show that the arc $(t^m, t^m, t^2)$ is in the closure of the
families (\ref{cA.sharc.thm}.2), provided we work in the space of formal arcs.
However, (\ref{short.defs.say}.1) is {\em not} a  deformation of short arcs
and  $(t^m, t^m, t^2)$ is a typical member of a new connected component
of $\sharc\bigl(0\in (xy=z^m)\bigr)$. 

By contrast, adding one more variable  kills this component.
For example, starting with the arc  $(t^m, t^m, t^2, 0)$
on $(xy=z^n)\subset \c^4$, we have
  deformations of short arcs
$$
(t,s)\mapsto  \bigl(t^i(t+s)^{m-i}, t^{m-i}(t+s)^{i}, t(t+s), ts\bigr).
\eqno{(\ref{short.defs.say}.2)}
$$
\end{say}

This example turns out to be typical and 
it is quite easy to modify the deformations in the proof of
Theorem \ref{cA.sharc.thm} to yield the following.

\begin{thm}\label{cA.sharc.sharc.thm}
 Let $X=(xy=f(z_1,\dots, z_n)\subset \c^{n+2}$ 
be a $cA$-type singularity.
 Assume that $\dim X\geq 3$ and $m:=\mult_0 f\geq 2$.

Then
 $\sharc(0\in X)$ has $(m-1)$
irreducible  components 
as in (\ref{cA.sharc.thm}.2).
\end{thm}

It is not always clear if a deformation $\Phi(t,s)$ 
is short or not. There is, however, one case when this is easy,
at least over a smaller disc $\dd_s(\epsilon)\subset \dd_s$.

\begin{lem} \label{easy.short.lem}
Let $\Phi(t,s)=\bigl(\Phi_1(t,s),\dots, \Phi_r(t,s)\bigr)$
be a deformation of arcs on $X\subset \c^r$.
Assume that $\Phi(t,0)$ is short and  $\Phi_i(t,s)$ is independent of $s$
and not identically zero 
for some $i$. Then  $\Phi(t,s_0):\bdd_t\to X$ is short for
$|s_0|\ll 1$.
\end{lem}

Proof. By assumption 
$\Phi(*,s_0)^{-1}(p)\subset \Phi_i(*,s_0)^{-1}(p)=\Phi_i(*,0)^{-1}(p)$
for every $s_0\in \dd_s$, thus there is a finite subset
$Z=\Phi_i(*,0)^{-1}(p)\subset \bdd_t$ such that
$$
\Phi^{-1}(0)\subset Z\times \dd_s\qtq{and} \Phi^{-1}(0)\cap (s=0)=\{(0,0)\}.
$$
Since $\Phi^{-1}(0) $ is closed, this implies that
$$
\Phi^{-1}(0)\cap \bigl(\bdd_t\times \dd_s(\epsilon)\bigr)
\subset \{0\}\times \dd_s(\epsilon)\qtq{for $0<\epsilon\ll 1$.}\qed
$$

\begin{say}[Proof of Theorem \ref{cA.sharc.sharc.thm}]
At the very beginning of the proof of Theorem \ref{cA.sharc.thm},
after a linear change of coordinates we may assume that
$z_1^m$ appears in $f$ with nonzero coefficient and
$\phi_2$ is not identically zero. Then the  construction
gives a deformation of short arcs by Lemma \ref{easy.short.lem}.

The deformations at the end of the proof were written to yield
short arcs. \qed
\end{say}

\section{A revised version of the Nash problem}\label{sec.revised}

As we saw, the   Nash map is not surjective in dimensions $\geq 3$.
In this section we  develop a revised version of the notion of
essential divisors. This leads to a smaller target for the
Nash map, so surjectivity should become more likely.
Our proposed variant of the Nash problem at least accounts for
all known counter examples.

We start with a reformulation of the original 
definition of essential divisors.

\begin{say} \label{reform.ess.say}
Let  $Y$ be  a complex variety and
$Z\subset Y$ a closed subset.  Let $\farc(Z\subset Y)$
denote the  scheme of formal arcs
$\phi:\spec \c[[t]]\to Y$ such that
$\phi(0)\in Z$. 

An easy but key observation is the following.
\medskip

\ref{reform.ess.say}.1. If $Y$ is smooth, then the irreducible components
of  $\farc(Z\subset Y)$ are in a natural one--to--one correspondence
with the irreducible components of $Z$. 
\medskip

We say that a
divisor $E$ over $Y$ is {\it essential} for $Z\subset Y$
if $E$ is obtained by blowing up one of the irreducible components of $Z$.
(For each irreducible component $Z_i\subset Z$, the blow-up
$B_ZY$ contains a unique divisor that dominates $Z_i$.)

The definition of essential divisors can now be reformulated as follows.
\medskip

\ref{reform.ess.say}.2.
Let $(x\in X)$ be a singularity. A divisor $E$  is
{\it essential} for $(x\in X)$ if $E$ is essential for
$\bigl(\supp\pi^{-1}(x)\subset Y\bigr)$
for every resolution $\pi:Y\to X$.
\end{say}

In order to refine the Nash problem, we need to understand
singular spaces for which the analog of (\ref{reform.ess.say}.1) still holds.

\begin{defn}[Sideways deformations]
Let $X$ be a variety (or an analytic space) and
$\phi:\spec \c[[t]]\to X$ a formal arc such that 
$\phi(0)\in \sing X$. 
A {\it sideways deformation} of  $\phi$ is a
morphism $\Phi:\spec \c[[t,s]]\to X$
such that 
$$
\Phi^* I_{\sing X}\subset (t,s)^m \qtq{for some $m\geq 1$}
$$
where $I_{\sing X}\subset \o_X$ is the ideal sheaf defining
$\sing X$. 

If $\Phi$ comes from a convergent arc
$\Phi^{\rm an}:\dd_t\times \dd_s\to X$
then this is equivalent to assuming that 
 for every $0\neq |s_0|\ll 1$ the nearby arc
$\Phi^{\rm an}(t,s_0)$ maps $\dd_t(\epsilon)$ to $X\setminus \sing X$
for some $0<\epsilon\leq 1$.

We say that $(x\in X)$ is {\it arc-wise Nash-trivial}
if every general arc in $\farc(x\in X)$
has a sideways deformation.
(By \cite{2012arXiv1201.6310F}, this implies that every  arc in $\farc(x\in X)$
has a sideways deformation.)
\end{defn}

\begin{comm} If $(x\in X)$ is  an isolated singularity with a
small resolution $\pi:X'\to X$ then every arc has a sideways deformation.
We can lift the arc to $X'$ and there move it away from the 
$\pi$-exceptional set. This is not very interesting and the
notion of essential divisors captures this phenomenon.

To exclude these cases, we are mainly interested in
arc-wise Nash-trivial singularities that do not have small
modifications.  
If arc-wise Nash-trivial singularities are log terminal
then assuming $\q$-factoriality captures this restriction, but in general
one needs to be careful of the difference between $\q$-factoriality
and having no small modifications.  

Also, in the few examples we know of, 
general arcs of every irreducible component
of $\farc(x\in X)$ have sideways deformations. If there are singularities
where sideways deformations exist only for some of the
irreducible components, the following outline needs to be suitably
modified.
\end{comm}

The main observation is that, for the purposes of the Nash problem,
 $\q$-factorial 
arc-wise Nash-trivial singularities should be considered as good as smooth
points.  The first evidence is the following straightforward analog of
(\ref{reform.ess.say}.1).

\begin{lem} Let $Y$ be  a complex variety with isolated,
arc-wise Nash-trivial singularities. Let
$Z\subset Y$ a closed subset that is the support of an effective
 Cartier divisor.
Then the irreducible components
of $\farc(Z\subset Y)$ are in a natural one--to--one correspondence
with the irreducible components of $Z$. \qed
\end{lem}

If $Z$ has lower dimensional  irreducible components, the situation
seems more complicated, but, at least in dimension 3, the
following seems to be the right generalization of (\ref{reform.ess.say}.1).

\begin{conj}  \label{easy.ess.conj}
Let $Y$ be  a 3--dimensional complex variety with isolated,
$\q$-factorial, 
arc-wise Nash-trivial singularities. Let
$Z\subset Y$ be a closed subset. 
Then the irreducible components
of $\farc(Z\subset Y)$ are in a natural one--to--one correspondence
with the union of the following two sets.
\begin{enumerate}
\item Irreducible components of $Z$.
\item Irreducible components of $\farc(p\in Y)$, where
$p\in Y$ is any
 singular point  such that $p\in Z$ and 
$\dim_pZ\leq 1$.
\end{enumerate}
\end{conj}

\begin{defn} Assumptions as in (\ref{easy.ess.conj}). 
A divisor over $Y$ is {\it essential} for $Z\subset Y$
if it corresponds to one of the irreducible components
of $\farc(Z\subset Y)$, as enumerated in (\ref{easy.ess.conj}.1--2).
\end{defn}

\begin{defn} Let $(x\in X)$ be a 3--dimensional, normal  singularity. 
A divisor $E$ over $X$ is called
{\it very essential}  for $(x\in X)$ if 
 $E$ is essential for
$\bigl(\supp\pi^{-1}(x)\subset Y\bigr)$
for every proper birational morphism $\pi:Y\to X$
where $Y$ has only   isolated,
$\q$-factorial, 
arc-wise Nash-trivial singularities.
(As in (\ref{first.ess.defn}),
 it is better to allow $Y$ to be an algebraic space.)
\end{defn}

It is easy to see that the Nash map
is an injection from the
irreducible components of $\farc(x\in X)$ into the set of
very essential divisors. 
One can  hope that there are no other obstructions.

\begin{prob}[Revised Nash problem]\label{reform.nash}
Is the   Nash map
 surjective  onto  the set of
very essential divisors for  normal 3-fold singularities? 
\end{prob}

As a first step, one should consider the following.

\begin{prob} In dimension 3, classify all 
$\q$-factorial, arc-wise Nash-trivial singularities.
\end{prob}

Hopefully they are all terminal and  a complete enumeration is
possible.  The papers \cite{MR2145316, MR2145317} contain 
  several results about partial resolutions of terminal singularities.

We treat two easy cases in
(\ref{sideways.cA}--\ref{sideways.quot}).
A positive solution of Question \ref{cDV.ques}
would imply that all isolated, 3-dimensional cDV singularities are
arc-wise Nash-trivial.

\begin{thm} \label{sideways.cA}
Let $(0\in X)$ be a $cA$-type singularity such that 
  $\dim \sing X\leq \dim X-3$. Then all 
 general arcs as in (\ref{cA.sharc.thm}.2) have sideways deformations.
\end{thm}

Proof. We use the notation of the proof of Theorem \ref{cA.sharc.thm}.

Since $\mult f\bigl(\phi_1(t), \dots, \phi_n(t)\bigr)=m$,
we see that $\mult \phi_j(t)=1$ for at least one index $j$.
We may assume that $j=1$ and $\phi_1(t)=t$.
 Thus, after the coordinate change 
$z_i\mapsto z_i-\phi_i(z_1)$  for $i=2,\dots,n$
and an additional general linear  coordinate change
among the $z_2,\dots, z_n$
we may assume that
\begin{enumerate}
\item $\phi_1(t)=t$,
\item $\phi_j(t)\equiv 0$ for $j>1$,
\item $\bigl(xy=g(z_1, z_2)\bigr)\subset \c^4$ has an isolated
singularity at the origin 
and $g(z_1,z_2)$ is  divisible neither by $z_1$ nor by $z_2$ 
where $g(z_1, z_2)=f(z_1, z_2, 0,\dots, 0)$. 
\end{enumerate}
 By  Lemma \ref{newton.lem}
there is an $r\geq 1$ such that
$$
g(t, s^r)=u(t,s)\prod_{i=1}^m \bigl(t-\sigma_i(s)\bigr).
$$
Since  $g(z_1,z_2)$ is not   divisible  by $z_1$
none of the $\sigma_i$ are identically zero.
Since $g(t, s) $ has an isolated singuarity at the origin and
is not  divisible  by $s$, 
$g(t, s^s) $ has an isolated singuarity at the origin.
Thus all the $\sigma_i(s)$ are distinct.

As before, for $j=1,2$ write  $\psi_j(t)=t^{a_j}v_j(t)$
where $v_j(0)\neq 0$. Note that $a_1+a_2=m$
and $u(t,0)=v_1(t)v_2(t)$.

Divide $\{1,\dots, m\}$ into two disjoint subsets $A_1, A_2$
such that $|A_j|=a_j$.
Finally set
$$
\Psi_1(t,s)=v_1(t)\cdot \prod_{i\in A_1} \bigl(t-\sigma_i(s)\bigr)
\qtq{and}
\Psi_2(t,s)=\frac{u(t,s)}{v_1(t)}
\cdot \prod_{i\in A_2} \bigl(t-\sigma_i(s)\bigr).
$$
Then
$$
\bigl(\Psi_1(t,s), \Psi_2(t,s), t, s^r, 0,\dots, 0\bigr)
$$
is a sideways deformation of 
$\bigl(\psi_1(t), \psi_2(t),t, 0,\dots, 0\bigr)$. \qed

\medskip

The opposite  happens for quotient singularities.

\begin{prop} \label{sideways.quot}
Let  $(0\in X):=\c^n/G$ be an isolated quotient singularity.
Then arcs with a sideways deformation are nowhere dense in
$\farc(0\in X)$. 
\end{prop}

Proof. 
Let   $\Phi:\spec \c[[t,s]]\to X$ be a sideways deformation
of an arc $\phi(t)=\Phi(t,0)$. By the purity of branch loci,
$\Phi$ lifts to an arc $\tilde \Phi:\spec \c[[t,s]]\to \c^n$.
In particular, $\phi:\spec \c[[t]]\to X$
lifts to $\tilde \phi:\spec \c[[t]]\to \c^n$.

By \cite{k-short}, such arcs constitute a connected component 
 of $\sharc(0\in X)$. 
We claim, however, that these arcs do not cover a whole
irreducible component of $\farc(0\in X)$. 

It is enough to show the
latter on some intermediate cover of $X$. The simplest is to use
 $(0\in Y):=\c^n/C$ where $C\subset G$ is any cyclic subgroup.

Set $r:=|C|$, fix a generator $g\in C$ and diagonalize its action as
$$
(x_1,\dots, x_n)\mapsto \bigl(\epsilon^{a_1}x_1,\dots, \epsilon^{a_n}x_n\bigr)
$$
where $\epsilon $ is a primitive $r$th root of unity.
%Set ${\mathbf v}:=\bigl(a_1/r, \dots, a_n/r\bigr)$.
Thus $Y$ is the toric variety corresponding to the 
free abelian group 
$$
N=\z^n+\z \bigl(a_1/r, \dots, a_n/r\bigr) \qtq{and the
cone} \Delta=\bigl(\q_{\geq 0}\bigr)^n.
$$
The Nash conjecture is true for toric singularities 
and by \cite[Sec.3]{MR2030097} the essential divisors are all toric
and  
correspond to interior vectors of
$N\cap \Delta$ that can not be written as the sum of
an interior vector of
$N\cap \Delta$ and of a nonzero vector of
$N\cap \Delta$.  In our case, all such vectors are
of the form 
$$
\bigl(\overline{ca_1}/r, \dots, \overline{ca_n}/r\bigr)
\qtq{for} c=1,\dots, r-1
$$
where  $\overline{ca_i} $ denotes remainder mod $r$.

Arcs that lift to $\c^n$ correspond to the vector
$(1,\dots, 1)$, which is not minimal. In fact
$$
(1,\dots, 1)=
\bigl(\overline{a_1}/r, \dots, \overline{a_n}/r\bigr)
+
\bigl(\overline{(r-1)a_1}/r, \dots, \overline{(r-1)a_n}/r\bigr). \qed
$$

%\bibliography{refs-main/refs}

\def\cprime{$'$} \def\cprime{$'$} \def\cprime{$'$} \def\cprime{$'$}
  \def\cprime{$'$} \def\cprime{$'$} \def\dbar{\leavevmode\hbox to
  0pt{\hskip.2ex \accent"16\hss}d} \def\cprime{$'$} \def\cprime{$'$}
  \def\polhk#1{\setbox0=\hbox{#1}{\ooalign{\hidewidth
  \lower1.5ex\hbox{`}\hidewidth\crcr\unhbox0}}} \def\cprime{$'$}
  \def\cprime{$'$} \def\cprime{$'$} \def\cprime{$'$}
  \def\polhk#1{\setbox0=\hbox{#1}{\ooalign{\hidewidth
  \lower1.5ex\hbox{`}\hidewidth\crcr\unhbox0}}} \def\cdprime{$''$}
  \def\cprime{$'$} \def\cprime{$'$} \def\cprime{$'$} \def\cprime{$'$}
\providecommand{\bysame}{\leavevmode\hbox to3em{\hrulefill}\thinspace}
\providecommand{\MR}{\relax\ifhmode\unskip\space\fi MR }
\providecommand{\MRhref}[2]{%
  \href{http://www.ams.org/mathscinet-getitem?mr=#1}{#2}
}
\providecommand{\href}[2]{#2}

\vskip1cm

\noindent Princeton University, Princeton NJ 08544-1000

{\begin{verbatim}jmjohnso@math.princeton.edu\end{verbatim}}

{\begin{verbatim}kollar@math.princeton.edu\end{verbatim}}

\end{document}